\documentclass[twoside]{ima}
\usepackage{amsmath}
\usepackage{amsfonts}
\usepackage{amssymb}
\usepackage{epsfig}
\usepackage{epstopdf}
\usepackage{url}

\newcommand{\CC}{\mathbb C}
\newcommand{\RR}{\mathbb R}
\newcommand{\QQ}{\mathbb Q}

\newcommand{\NN}{\mathbb N}
\newcommand{\PP}{\mathbb P}

\begin{document}

\markboth{BERND STURMFELS}{ALGEBRAIC STATISTICS}

\title{OPEN PROBLEMS IN ALGEBRAIC STATISTICS}

\author{BERND STURMFELS\thanks{University of California, Berkeley, CA 94720, USA, {\tt bernd@math.berkeley.edu}}
}
 \maketitle

\date{\today}

\begin{abstract}      
Algebraic statistics is concerned with the study of probabilistic models
and techniques for statistical inference using methods from algebra
and geometry. This article presents a list of open mathematical problems
in this emerging field, with main emphasis on
graphical models with hidden variables, maximum likelihood estimation,
and multivariate Gaussian distributions.
These are notes from a lecture presented at the IMA in
Minneapolis during the 2006/07 program on
Applications of Algebraic Geometry. 

\end{abstract} 


{\AMSMOS 13P10, 14Q15, 62H17, 65C60
\endAMSMOS}

 \section{Introduction}
 
 This article is based on a lecture given in March 2007
 at the workshop on {\em Statistics, Biology and Dynamics} held
 at the Institute for Mathematics and its Applications (IMA) in Minneapolis
 as part of the  2006/07 program on {\em Applications of Algebraic Geometry}.
 In four sections we present mathematical problems whose solutions would likely become
 important contributions to the emerging interactions between algebraic geometry
 and computational statistics. Each of the four sections starts out with a ``specific problem''
which plays the role of representing the broader research agenda. The latter is summarized in a
``general problem''.

Algebraic statistics is concerned with the study of probabilistic models
and techniques for statistical inference using methods from algebra
and geometry. The term was coined in the book of 
Pistone, Riccomagno and Wynn \cite{PRW}
and subsequently developed for biological applications in \cite{ASCB}.
Readers from statistics will enjoy the introduction and review recently given by
Drton and Sullivant \cite{DS}, while readers from algebra will find
various points of entry
 cited in our discussion and listed among our references. 
 
 \section{Graphical Models with Hidden Variables}

Our first question concerns three-dimensional contingency tables
$(p_{ijk})$ whose indices $i,j,k$ range over a set of four elements,
such as the set $ \{{\tt A},{\tt C},{\tt G},{\tt T}\}$ of DNA bases.

\smallskip 

{\bf Specific Problem: }
{\em Consider the 
variety of $4 {\times} 4 {\times} 4$-tables
of tensor rank at most $4$.
There are certain known polynomials
of degree at most nine which vanish on this variety.
Do they suffice to cut out the variety?}

\smallskip

\noindent
This particular open problem appears in \cite[Conjecture 3.24]{ASCB},
and it here serves as a placeholder for the following broader direction of inquiry.

\smallskip

{\bf General Problem:}
{\em Study the geometry and commutative algebra
of graphical models with hidden random variables.
Construct these varieties by gluing  familiar 
secant varieties, and by applying representation theory.}

\smallskip

We are interested in statistical models for discrete data 
which can be represented by polynomial constraints. As is
customary in algebraic geometry, we consider varieties
over the field of complex numbers, with the tacit
understanding that statisticians mostly care about
points whose coordinates are real and non-negative.
The model referred to in the Specific Problem
lives in the $64$-dimensional space
 $\,\, \CC^4 \otimes \CC^4 \otimes \CC^4$ 
of $4 {\times} 4 {\times 4}$-tables $(p_{ijk})$,
where $i,j,k  \in \{{\tt A},{\tt C},{\tt G},{\tt T}\}$. 
It has the parametric representation
\begin{equation}
\label{para444}
 \begin{matrix}
p_{ijk} \quad = \quad &
 \rho_{{\tt A}i} \cdot  \sigma_{{\tt A}j} \cdot \theta_{{\tt A}k} \, + 
 \rho_{{\tt C}i} \cdot  \sigma_{{\tt C}j} \cdot \theta_{{\tt C}k} +  \\ &
 \rho_{{\tt G}i} \cdot  \sigma_{{\tt G}j} \cdot \theta_{{\tt G}k} \,+ \,
 \rho_{{\tt T}i} \cdot  \sigma_{{\tt T}j} \cdot \theta_{{\tt T}k}. \,\,
\end{matrix}
\end{equation}
Our problem is to compute the
homogeneous prime ideal $I$ of 
all polynomials which vanish on this model.
The desired ideal $I$ lives in the polynomial ring 
$\,\QQ \bigl[ p_{{\tt A} {\tt A} {\tt A}},
 p_{{\tt A} {\tt A} {\tt C}}, 
 p_{{\tt A} {\tt A} {\tt T}}, 
\ldots, p_{{\tt T} {\tt T} {\tt G}},
 p_{{\tt T} {\tt T} {\tt T}} \bigr] $
with $64$ unknowns. In principle, one can compute
generators of $I$ by applying Gr\"obner bases methods
to the parametrization (\ref{para444}).
However, our problem has $64$ probabilities and
$48$ parameters, and it is simply too big
for the kind of computations which were performed
in \cite[\S 3.2]{ASCB} using
 the software package {\tt Singular} \cite{singular}. 

Given that Gr\"obner basis methods appear to be
too slow for any problem size which is actually
relevant for real data, skeptics
may wonder why a statistician should bother 
learning the language of ideals and varieties. 
One possible response to the practitioner's legitimate question
{\em ``Why (pure) mathematics?''} is offered by the following quote due to
Henri Poincar\'e:

\smallskip

\centerline{\em
``Mathematics is the Art of Giving the 
  Same Name to Different Things''.}

 \medskip
 
Indeed, our prime ideal $I$ gives the same name to
the following things:
\begin{itemize}
\item the set of $4 {\times} 4 {\times} 4$-tables of tensor rank $\leq 4$,
\item  the mixture of four models for three independent random variables,
\item the naive Bayes model with four classes,
\item the conditional independence model
$\,[X_1 \perp \!\!\! \perp  X_2 \perp \!\!\! \perp X_3 \,|\, Y \,]$,
\item the fourth secant variety of the Segre variety $\PP^3 {\times} \PP^3 {\times} \PP^3 $,
\item the general Markov model for the phylogenetic tree $K_{1,3}$,
\item superposition of four pure states in a quantum system \cite{Brody, Hey}.
\end{itemize}
These different terms have been used in the literature
for the geometric object represented by  (\ref{para444}).
The concise language of commutative algebra
and algebraic geometry can be an effective
channel of communication for the different communities
of statisticians, computer scientists, physicists, engineers
and biologists, all of whom have encountered formulas like (\ref{para444}).

The generators of lowest degree in our ideal $I$ have degree five, 
and the known generators of highest degree have degree nine.
The analysis of Landsberg and Manivel in \cite[Proposition 6.3]{LM} 
on $3 {\times} 3 {\times} 4$-tables
of tensor rank four implies the existence
of additional ideal generators of degree six in $I$. This analysis had 
been overlooked by the authors of \cite{ASCB} when they formulated
their Conjecture 3.24. Readers of \cite[Chapter 3]{ASCB} are herewith kindly asked to replace
{\em ``of degree $5$ and $9$''} by {\em ``of degree at most $9$''}.

In what follows we present the known minimal generators
of degree five and nine in our prime ideal $I$, and we postpone a more detailed 
discussion of the Landsberg-Manivel sextics 
in \cite[Proposition 6.3]{LM} to a future study.

Consider any $3 \times 4 \times 4$-subtable $(p_{ijk})$
and let $A,B,C$ be the $4 {\times} 4$-slices gotten by fixing $i$.
To be precise, the entry of the $4 \times 4$-matrix $A$  in row $j$ and column $k$
equals $p_{{\bf A}jk}$, 
the entry of $B$ in row $j$ and column $k$  equals $p_{{\bf C}jk}$, and 
the entry of $C$ in row $j$ and column $k$  equals $p_{{\bf G}jk}$.
We can check
that the following identity of $4 {\times} 4$-matrices holds
for all tables in our model, provided the matrix $B$ is invertible:
$$ A \cdot B^{-1} \cdot C \,\, = \,\, C \cdot B^{-1} \cdot A $$
After clearing the denominator $\,{\rm det}( B)$, we can write this
identity as
\begin{equation}
\label{adjointeqn}
 A \cdot {\rm adj}(B) \cdot C \, - \,C \cdot {\rm adj}(B) \cdot A \quad = \quad 0, 
 \end{equation}
where ${\rm adj}(B) = {\rm det}(B) \cdot B^{-1}$ is the adjoint matrix of $B$.
The matrix entries on the left hand side give $16$ quintic polynomials 
which lie in our prime ideal $I$.  Each matrix entry is a polynomial
with $180$ terms which involve only $30$ of the $64$ unknowns.
 For example, the upper left entry looks like this:
\begin{eqnarray*}
& \phantom{-} p_{\bf AAC} p_{\bf CCA} p_{\bf CGG} p_{\bf CTT} p_{\bf GAA}
-p_{\bf AAC} p_{\bf CCA} p_{\bf CGT} p_{\bf CTG} p_{\bf GAA} \\
& -p_{\bf AAC} p_{\bf CCG} p_{\bf CGA} p_{\bf CTT} p_{\bf GAA}
+ p_{\bf AAC} p_{\bf CCT} p_{\bf CGA} p_{\bf CTG} p_{\bf GAA} \\
& \, + \quad  \cdots \cdots \,\hbox{(175 terms)} \cdots \cdots \quad
-\,\, p_{\bf ATA} p_{\bf CAG} p_{\bf CCC} p_{\bf CGA} p_{\bf GAT}.
\end{eqnarray*}

We note that there are no non-zero polynomials of degree $\leq 4$ in the ideal $I$.
This follows from general results on secant varieties \cite{Cat, LM1}.

An explicit linear algebra computation reveals that all polynomials of degree five in $I$
are gotten from the above construction by relabeling and considering all
subtables of format $3 {\times} 4 {\times} 4$, format
$4 {\times} 3 {\times} 4$ and format $4 {\times} 4 {\times} 3$,
and by applying the natural action of the group
$\,GL(\CC^4) \times GL(\CC^4) \times GL(\CC^4)$
on $4 {\times} 4 {\times} 4$-tables. This action leaves the ideal $I$ fixed.
We identify the representation of this group on the
space of quintics in $I$.

\smallskip

\begin{proposition} \label{propdeg5}
The space of quintic polynomials in the prime ideal $I$ of (\ref{para444})
has dimension $1728$. As a $GL(\CC^4)^3 $-module, it is isomorphic to
$$ 
\begin{matrix} \quad \,\,\,
S_{311} (\CC^4) \otimes S_{2111}(\CC^4) \otimes S_{2111}(\CC^4)  \\ \oplus \,\,\,\,
S_{2111} (\CC^4) \otimes S_{311}(\CC^4) \otimes S_{2111}(\CC^4)  \\ \, \oplus \,\,\,\,
S_{2111} (\CC^4) \otimes S_{2111}(\CC^4) \otimes S_{311}(\CC^4)  .
\end{matrix}
$$
\end{proposition}

Here $S_\lambda(\CC^4)$ denotes the {\em Schur modules}
which are the irreducible representations of $GL(\CC^4)$.
We refer to \cite{FH} for the relevant basics on representation theory
of the general linear group,
and to \cite{LM1, LM2, LW} for more detailed information about
the specific modules under consideration here.

The known invariants of degree nine are also obtained by a similar construction.
Consider any $3 \times 3 \times 3$-subtable  $(p_{ijk})$
and denote the three slices of that table by $A$, $B$ and $C$.
We now consider the  $3 {\times} 3$-determinant
\begin{equation}
\label{strassenmatrix}
 {\rm det} ( A \cdot B^{-1} \cdot C \, - \, C \cdot B^{-1} \cdot A) .
\end{equation}
The denominator of the rational function (\ref{strassenmatrix}) is ${\rm det}(B)^2$
and not ${\rm det}(B)^3$ as one might think on first glance.
The numerator of (\ref{strassenmatrix}) is a homogeneous polynomial
of degree nine with $9216$ terms which remains invariant under permuting
$A$, $B$ and $C$. 
This homogeneous polynomial of degree nine
lies in the ideal $\,I\,$ and is known as the {\em Strassen invariant}.

\smallskip

\begin{proposition} \label{propdeg9}
The  $GL(\CC^4)^3 $-submodule of the
degree $9$ component $I_9$  generated by
the Strassen invariant  is not contained in
the ideal $\langle I_5 \rangle$ generated by the quintics in Proposition \ref{propdeg5}.
This module  has vector space dimension $8000$ and it is isomorphic to
the representation
$$ S_{333} (\CC^4) \otimes S_{333}(\CC^4) \otimes S_{333}(\CC^4)  .$$
\end{proposition}

The first appearance of the Strassen invariant in algebraic statistics
was \cite[Proposition 22]{GSS}.
A conceptual study of the matrix construction 
$\, A  B^{-1}  C \, - \, C B^{-1}  A \,$
was undertaken by Landsberg and Manivel in \cite{LM2}.

The Specific Problem at the beginning of this section
plays a pivotal role also in algebraic phylogenetics \cite{salmon, AR1, AR2}.
Our model (\ref{para444}) is known there as the
general Markov model on a tree with three leaves
branching off directly from the root.
Allman and Rhodes \cite[\S 6]{AR1} showed that
phylogenetic invariants which cut out the general Markov model on any
larger binary rooted tree can be constructed
from the generators of our ideal $I$ by a gluing process.
The invariants of degree five and nine arising from 
(\ref{adjointeqn}) and (\ref{strassenmatrix}) are
therefore basic building blocks for phylogenetic
invariants on arbitrary trees whose nodes
are labeled with the four letters ${\bf A}$, ${\bf C}$, ${\bf G}$~and~${\bf T}$.

In her lecture at the same IMA conference in March 2007,
Elizabeth Allman \cite{salmon} offered an extremely attractive prize
for the resolution of the Specific Problem.
She offered to personally catch and smoke wild salmon 
from the Copper River, located in her ``backyard'' in Alaska,
and ship it to anyone who will determine a minimal
generating set of the prime ideal $I$.

\smallskip

In Propositions \ref{propdeg5} and \ref{propdeg9},
we emphasized the language of
representation theory in characterizing the
defining equations of graphical statistical models.
This methodology is a main focus in the forthcoming
book by J.M.~Landsberg and Jason Morton, which advocates
the idea of using Schur modules $S_\lambda(\CC^n)$
in the description of such models. Morton's key insight
 is that this naturally generalizes  conditional independence, 
 the current language of choice for characterizing 
 graphical models. Conditional independence statements
 can be interpreted as a convenient shorthand
for large systems of quadratic equations;
see~\cite[\S 4.1]{GMS} or \cite[Proposition 8.1]{Solving}.

In the absence of hidden random variables, the quadratic equations
expressed implicitly by conditional independence are
sufficient to characterize graphical models.
This is the content of the {\em Hammersley-Clifford Theorem} 
(see e.g.~\cite[Theorem 4.1]{GMS} or \cite[Theorems 1.30 and 1.33]{ASCB}).
However, when some of the random variables
in a graphical model are hidden then the situation
becomes much more complicated. We believe
that representation theory of the general linear group
can greatly enhance the conditional independence calculus
which is so widely used by graphical models experts.
The representation-theoretic notation was here illustrated for a tiny graphical model, 
having three observed random variables
and one hidden random variable, all four having the same state space
$\{{\bf A}, {\bf C}, {\bf G}, {\bf T} \}$.

\medskip

\section{Maximum Likelihood Estimation} 

In this section we discuss topics concerning the algebraic approach to
 maximum likelihood estimation \cite[\S 3.3]{ASCB}.
 The following open problem was published in \cite[Problem 13]{HKS}.

\smallskip

\noindent
{\bf Specific Problem: } {\em
Find a geometric characterization of those projective varieties
whose maximum likelihood degree (ML degree) is equal to one.}
\smallskip

This question and others raised in \cite{CHKS, HKS}
are just the tip of an iceberg:

\smallskip

\noindent
{\bf General Problem:} {\em
Study the geometry of maximum likelihood estimation
for algebraic statistical models. }

\smallskip

Here algebraic statistical models are regarded as projective varieties.
A model has ML degree one if and only if
its maximum likelihood estimator is a rational function
of the data. Models which have this property tend to be
very nice. For instance, in the special context of undirected
graphical models (Markov random fields), the property
of having ML degree one is equivalent to the statement
that the graph is decomposable \cite[Theorem 4.4]{GMS}.
For toric varieties, our question was featured in \cite[Problem 8.23]{Solving}.

It is hoped that the ML degree is related to
convergence properties of numerical algorithms used by statisticians,
such as iterative proportional scaling or the EM algorithm,
but no systematic study in this direction has yet been undertaken.
In general, we wish to learn how statistical features
of a model relate to geometric properties of the
corresponding variety.

Here are the relevant definitions for our problems.
We fix the complex projective space $\PP^n$ with coordinates
$(p_0:p_1:\cdots: p_n)$.  The coordinate
$p_i$ represents the probability of the $i$th event.
The $n$-dimensional probability simplex
is identified with the set $\,\PP^n_{\geq 0}\,$
of points in $\PP^n$ which have
non-negative real coordinates. The data comes in the  form of 
a non-negative integer vector $ (u_0,u_1,\ldots,u_n) \in \NN^{n+1}$.
Here $u_i$ is the number of times the $i$th event was observed.
The corresponding {\em likelihood function} is defined as
 \begin{equation}
 \label{likelihood}
L(p_0,p_1,\ldots,p_n) \,\,\,\, = \,\,\,\,
\frac{ {p_0}^{u_0}
 \cdot {p_1}^{u_1}
  \cdot {p_2}^{u_2}  \,\cdot \cdots \cdot \,
  {p_n}^{u_n}} {
 (p_0 {+} p_1 {+}        \cdots {+} p_n)^{u_0+u_1+ \cdots +u_n}}.
 \qquad
 \end{equation}
 Statistical computations are typically done in 
  affine $n$-space specified by $p_0 + p_1  + \cdots + p_n = 1$,
  where the denominator of $L$ can be ignored. However,
the denominator is needed in order for $L$ to be a
well-defined rational function on $ \PP^n$. 
The unique critical point of the likelihood function $L$ is at $(u_0: u_1: \cdots : u_n)$,
and this point is the global maximum of $L$ over $\PP^n_{\geq 0}$.
By a {\em critical point} we mean any point at which the gradient of $L$
vanishes.

An {\em algebraic statistical model} is represented by a subvariety $\mathcal{M}$
of the projective space $\PP^n$. 
The model itself is the intersection of $\mathcal{M}$
with the probability simplex $\PP^n_{\geq 0}$.
The {\em ML degree} of the variety $\mathcal{M}$ is
the number of complex critical points of the restriction
of the likelihood function $L$ to $\mathcal{M}$.
Here we disregard singular points of $\mathcal{M}$, we 
only count critical points that are not
poles or zeros of $L$, and $ u_0,u_1,\ldots,u_n$ are assumed to be generic.
If $\mathcal{M}$ is smooth and the 
divisor on $\mathcal{M}$ defined by $L$ has normal crossings then 
there is a geometric characterization of the ML degree,
derived in the paper  \cite{CHKS} with Catanese, Ho\c{s}ten and Khetan.
The assumptions of smoothness and normal crossing are
very restrictive and almost never satisfied for models of
statistical interest. In general, to understand the ML degree
will require invoking some resolution of singularities
and its algebraic underpinnings.

We illustrate the computation of the ML degree for the
case when $\mathcal{M}$ is a plane curve.
Here $n=2$ and $\mathcal{M}$ is the zero set of
a homogeneous polynomial  $\, F(p_0,p_1,p_2) $.
Using  Lagrange multipliers or \cite[Proposition 2]{HKS},
we derive that the condition
for $(p_0:p_1: p_2)$ to be a critical point of
the restriction of $L$ to $\mathcal{M}$ is equivalent to
the system of two equations
$$ F(p_0,p_1,p_2) \quad = \quad
{\rm det}
\begin{pmatrix}
\,\,u_0 \, & p_0 \,& p_0 \cdot \partial{F}/\partial{p_0} \\
\,\,u_1 \, & p_1 \,& p_1 \cdot \partial{F}/\partial{p_1} \\
\,\,u_2 \, & p_2 \,& p_2 \cdot \partial{F}/\partial{p_2} 
\end{pmatrix} \quad = \quad 0.
$$
For a general polynomial $F$ of degree $d$, these equations
will have $d(d+1)$ solutions, by B\'ezout's Theorem.
Moreover, all of these solutions satisfy
\begin{equation}
\label{genericity} p_0 \cdot p_1 \cdots p_n \cdot (p_0 + p_1 + \cdots + p_n) \,\, \not= \,\, 0,
\end{equation}
and we conclude that the ML degree of a general plane curve
of degree $d$ is equal to $d(d+1)$. However, that number can drop considerably
for special curves. For instance, while the
ML degree of a general plane quadric equals six, the special
quadric $\, \{p_1^2 = \lambda p_0 p_2 \}\,$
has ML degree two for $\lambda \not= 4$, and it
has ML degree one for $\lambda = 4$.
Thus, returning to the Special Problem,
our first example of a variety of ML degree one is the plane curve defined by
\begin{equation}
\label{2x2symm}
 F \quad = \quad
{\rm det} \begin{pmatrix} 
2 p_0 & p_1 \\
   p_1 & 2 p_2 
   \end{pmatrix}   .
   \end{equation}
Biologists know this as the  {\em Hardy-Weinberg curve}, with the
 parametrization
\begin{equation}
\label{HW} p_0 \,=\, \theta^2 \,,\quad
   p_1\, =\, 2 \theta(1-\theta)\, , \quad
   p_2  \,=\, (1-\theta)^2.   
   \end{equation}
 The unique critical point of the likelihood function $L$ on this curve equals
$$ \bigl(   \,  (2 u_0+u_1)^2\,:\, 2(2 u_0 {+} u_1) (u_1 {+} 2 u_2) \, : \,  (u_1+2u_2)^2 \bigr) . $$

Determinantal varieties arise naturally in statistics.
They are the models $\mathcal{M}$ that are specified by
imposing  rank conditions on a matrix of unknowns.
A first example is the model (\ref{HW})
for two i.i.d.~binary random variables.
For a second example we consider the general $3 \times 3$-matrix
\begin{equation}
\label{3x3matrix} P \quad = \quad  \begin{pmatrix}
p_{00} & p_{01} & p_{02} \\
p_{10} & p_{11} & p_{12} \\
p_{20} & p_{21} & p_{22} 
\end{pmatrix} 
\end{equation}
which represents two ternary random variables.
The independence model for these two random variables
is the variety of rank one matrices. This model also has
ML degree one, i.e., the maximum likelihood estimator 
is a rational function in the data. It is given by the $3 \times 3$-matrix
whose entry in row $i$ and column $j$ equals
$\, (u_{i0} + u_{i1} + u_{i2}) \cdot (u_{0j} + u_{1j} + u_{2j}) $.

By contrast, consider the {\em mixture model}
based on  two ternary random variables. It consists
of all matrices $P$ of rank at most two. Thus this model is
the hypersurface defined  by the cubic polynomial
$\, F \, = \, {\rm det}(P)$. Explicit computation shows that the 
ML degree of this hypersurface is ten.
In general, it remains an open problem to find a formula,
in terms of $m,n$ and $r$, for the
 ML degree of the variety of
$m {\times} n$-matrices of rank $\leq r$.

The first interesting case arises when $m=n=4$ and $r=2$.
At present we are unable to solve the likelihood equations 
for this case symbolically. The following concrete 
 biology example was proposed in \cite[Example 1.16]{ASCB}:

\smallskip

{\em 
``Our data are two aligned DNA sequences ...
  $$
 \begin{matrix}{\tt A T C A C C A A A C A T T G G G A T G C C T G T G
  C A T T T G C A A G C G G C T} \\{\tt A T G A G T C T T A A A C G C
  T G G C C A T G T G C C A T C T T A G A C A G C G}
\end{matrix}
$$
.. test the hypothesis that these two sequences were generated
by DiaNA using one biased coin and four tetrahedral dice....''}

\smallskip

Here the model $\mathcal{M}$ consists of all (positive)
 $4 {\times} 4$-matrices $ ( p_{ij})$ of rank at most two.
 In the given alignment, each match occurs four times
 and each mismatch occurs two times. Hence the likelihood function (\ref{likelihood}) equals
 $$ \!\!\!\! L \,\,\, = \,\,\,
(\prod_{i} p_{ii})^{4} \cdot
(\prod_{i \not= j} p_{ij})^{2} \cdot (\sum_{i,j} p_{ij})^{-40} .$$
Based on experiments with the EM algorithm, 
we conjectured that
the matrix $\,\bigl(\hat{p}_{ij}\bigr)  = 
\frac{1}{40} \begin{pmatrix}
3 & 3 & 2 & 2 \\
3 & 3 & 2 & 2 \\
2 & 2 & 3 & 3 \\
2 & 2 & 3 & 3 
\end{pmatrix}$ 
is a global maximum of the likelihood function $L$.
In the {\em Nachdiplomsverlesung} (postgraduate course) which I held at ETH Z\"urich in
the summer of 2005, I offered a cash
prize of 100 Swiss Francs for the resolution of this
very specific conjecture, and this prize remains unclaimed 
and is still available at this time (August 2007).

 The state of the art on
this {\em 100 Swiss Francs Conjecture}
is the work of Hersh which originated
in March 2007 at the IMA. She proved a 
range of constraints on the
maximum likelihood estimates of
determinantal models, especially
when the data $u_{ij}$ have symmetry.
A discussion of these ideas appears in
Hersh's paper with Fienberg, Rinaldo and Zhou \cite{HFRZ}.
That paper gives an exposition of MLE for
determinantal models aimed at statisticians.

\medskip

\section{Gaussian Conditional Independence Models}

The early literature on algebraic statistics,
including the book \cite{ASCB}, dealt primarily
with discrete random variables (binary, ternary,$\ldots$).
The set-up was as described in the previous two sections.
We now shift gears and consider multivariate Gaussian distributions.
For continuous random variables, we must work
in the space of model parameters in order to apply
algebraic geometry.
The following concrete problem concerns
Gaussian distributions on $\RR^5$.

\smallskip

\noindent {\bf Specific Problem:} {\em 
Which sets of almost-principal minors can be zero  
for a positive definite symmetric $5 {\times} 5$-matrix?}

\smallskip

The general question behind this asks for
characterization of all
  conditional independence models which can be
  realized by Gaussians on $\RR^n$.

\smallskip

\noindent {\bf  General Problem:} {\em
Study the geometry of conditional independence models
for multivariate Gaussian random variables.}

\smallskip

The state of the art on these problems appears in the work of
Franti\v{s}ek Mat\'u\v{s} and his collaborators. In particular,
Mat\'u\v{s}' recent paper with Ln\v{e}ni\v{c}ka \cite{LM} on
{\em representation of gaussoids} solves our Specific Problem for symmetric 
$4 {\times} 4$-matrices. Sullivant's construction in \cite{Sul} complements that work.
For more information see also the article by \v{S}ime\v{c}ek \cite{Sim}.

Let us begin, however, with some basic definitions. Our aim
is to discuss these problems in a self-contained manner.
A {\em multivariate Gaussian} distribution on $\RR^n$ with mean zero
is specified by its covariance matrix
$\,\Sigma = (\sigma_{ij})$. The $n {\times} n$-matrix $\Sigma$ is symmetric
and it is {\em positive definite}, which means that all its $2^n$ principal minors are
positive real numbers.

An {\em almost-principal minor} of $\,\Sigma\,$
is a subdeterminant which has row indices $ \{i\} \cup K \,$
and column indices $\,\{j\} \cup K \,$
for some $K \subset \{1,\ldots,n\}$ and $i,j \in \{1,\ldots,n\} \backslash K$.
We denote this subdeterminant by $[\,i \!\perp \!\!\! \perp\! j \, | K \,]$.
For example, if $n=5$, $i=2, j = 4$ and $K = \{1,5\}$ then
the corresponding almost-principal minor of the symmetric
$5 {\times} 5$-matrix $\Sigma$ equals
$$
[\,2 \!\perp \!\!\! \perp\! 4 \, | \{1,5\} \,] \quad = \quad
{\rm det}
\begin{pmatrix}
\sigma_{24} & \sigma_{12} & \sigma_{25} \\
\sigma_{14} & \sigma_{11} & \sigma_{15} \\
\sigma_{45} & \sigma_{15} & \sigma_{55} 
\end{pmatrix}
$$

Our notation for almost-principal minors is justified by their
intimate connection to conditional independence, 
expressed in the following lemma.
We note that the almost-principal minors
are referred to as {\em partial covariance} (or, if renormalized, {\em partial correlations})
 in the statistics literature.

\smallskip

\begin{lemma} \label{lem:covariance}
The subdeterminant $[\,i \!\perp \!\!\! \perp j\,| K \,] $ is zero for
a positive definite symmetric $n {\times} n$-matrix $\Sigma$ if and only if,
for the Gaussian random variable $X$ on $\RR^n$ with covariance matrix $\Sigma$,
the random variable
$ X_i  $ is independent of the random variable $X_j$  given
the joint variable $X_K$.
\end{lemma}

\smallskip

\begin{proof}
See \cite[Equation(5)]{DSS},  \cite[Section 1]{Mat}, or
\cite[Proposition 2.1]{Sul}.
\end{proof}

\medskip

Let  $\,{\tt PD}_n\,$ denote the $\binom{n+1}{2}$-dimensional
cone  of positive definite symmetric $n {\times} n$-matrices.
Note that this cone is open.
A {\em Gaussian conditional independence model}, or
{\em GCI model} for short,
is any semi-algebraic subset of the cone  $\,{\tt PD}_n\,$ 
which can be defined by polynomial equations of the form
\begin{equation}
\label{CIeqn}
 \,[\,i \!\perp \!\!\! \perp j\,| K \,] \quad = \quad 0 .
 \end{equation}
 In algebraic geometry, we simplify matters by
 studying the complex algebraic varieties defined
 by equations of the form (\ref{CIeqn}). Of course, what we are
 particularly interested in is the
 real locus of such a complexified GCI model, and how it intersects the 
 positive definite cone $\,{\tt PD}_n\,$ and its closure.

\smallskip

As an illustration of algebraic reasoning for Gaussian
conditional independence models, we examine
an example taken from \cite{Sul}.
Let $n=5$ and consider the GCI model given by the five quadratic polynomials
$$ \begin{matrix}
\,[\,1 \!\perp \!\!\! \perp 2 \,\,|\, \{3\} \,] && = && \sigma_{12} \sigma_{33} - \sigma_{13} \sigma_{23} \\
\,[\,2 \!\perp \!\!\! \perp 3 \,\,|\, \{4\} \,] && = && \sigma_{23} \sigma_{44} - \sigma_{24} \sigma_{34} \\
\,[\,3 \!\perp \!\!\! \perp 4 \,\,|\, \{5\} \,] && = && \sigma_{34} \sigma_{55} - \sigma_{35} \sigma_{45} \\
\,[\,4 \!\perp \!\!\! \perp 5 \,\,|\, \{1\} \,] && = && \sigma_{45} \sigma_{11} - \sigma_{14} \sigma_{15} \\
\,[\,5 \!\perp \!\!\! \perp 1 \,\,|\, \{2\} \,] && = && \sigma_{15} \sigma_{22} - \sigma_{25} \sigma_{12} 
\end{matrix} $$
This variety is a complete intersection (it has dimension ten) in the
$15$-dimensional space of symmetric $5 {\times} 5$-matrices.
Primary decomposition reveals that it is the union of precisely two
irreducible components, namely,
\begin{itemize}
\item  the linear space $\,\{ \,\sigma_{12} = 
 \sigma_{23} =  \sigma_{34} =  \sigma_{45} =  \sigma_{15} = 0 \,\}$, and
\item the toric variety defined by the five quadrics plus the extra equation
\begin{equation}
\label{magicquintic}
\sigma_{11} \sigma_{22} \sigma_{33} \sigma_{44} \sigma_{55} \,\,
= \,\, \sigma_{13} \sigma_{14} \sigma_{24} \sigma_{25} \sigma_{35}. 
\end{equation}
\end{itemize}
All matrices in the open cone $\,{\tt PD}_5\,$ satisfy
the inequalities $\,\sigma_{ii} >  0 \,$ and
$$ \sigma_{11} \sigma_{33} > \sigma_{13}^2 \, ,\,
 \sigma_{22} \sigma_{44} > \sigma_{24}^2 \,,\, 
 \sigma_{33} \sigma_{55} > \sigma_{35}^2 \,,\,
   \sigma_{44} \sigma_{11} > \sigma_{14}^2 \, ,\,
 \sigma_{55} \sigma_{22} > \sigma_{25}^2 . $$
 Multiplying the left hand sides and right hand sides respectively, we find
 $$ 
 \sigma_{11}^2 \sigma_{22}^2 \sigma_{33}^2 \sigma_{44}^2 \sigma_{55}^2 \, > \,
 \sigma_{13}^2 \sigma_{14}^2 \sigma_{24}^2 \sigma_{25}^2 \sigma_{35}^2 . 
 $$
This is a contradiction to the equation (\ref{magicquintic}), and 
we conclude that the intersection of our GCI model  with ${\tt PD}_5$
is contained in the linear space $\,\{ \,\sigma_{12} = 
 \sigma_{23} =  \sigma_{34} =  \sigma_{45} =  \sigma_{15} = 0 \,\}$.
The vanishing of the off-diagonal entry $\sigma_{ij}$
means that $X_i$ is independent of $X_j$, or, in symbols,
$\,[\,i \!\perp \!\!\! \perp \! j \,]$.
Our algebraic computation thus implies the following axiom for
GCI models.

\smallskip

\begin{corollary} \label{cororo}
Suppose the conditional independence statements 
$ \,[\,1 \!\perp \!\!\! \perp 2 \,\,|\, \{3\} \,] $, 
$\,[\,2 \!\perp \!\!\! \perp 3 \,\,|\, \{4\} \,] $,
$\,[\,3 \!\perp \!\!\! \perp 4 \,\,|\, \{5\} \,] $,
$\,[\,4 \!\perp \!\!\! \perp 5 \,\,|\, \{1\} \,] $,
$\,[\,5 \!\perp \!\!\! \perp 1 \,\,|\, \{2\} \,] $ hold
for some multivariate Gaussian distribution.
Then also the following five statements must hold:
$[\,1 \!\perp \!\!\! \perp \! 2 \,]$,
$[\,2 \!\perp \!\!\! \perp \! 3 \,]$,
$[\,3 \!\perp \!\!\! \perp \! 4 \,]$,
$[\,4 \!\perp \!\!\! \perp \! 5 \,]$
and $ [\,5 \!\perp \!\!\! \perp \! 1 \,]$.
\end{corollary}

\smallskip

Let us now return to the question
{\em ``which almost-principal minors can simultaneously vanish
for a positive definite symmetric $n \times n$-matrix?'' }
Corollary \ref{cororo} gives a necessary condition for $n = 5$.
We next discuss the answer to our question for $n \leq 4$.
For $n = 3$, the necessary and sufficient conditions
are given (up to relabeling) by the following four axioms:

\begin{itemize}
\item[(a)] $ \,[\,  1 \!\perp \!\!\! \perp  2  \,] \,$ and
$ \,[\,  1 \!\perp \!\!\! \perp 3  \,\,|\, \{ 2   \} \,] \,$ \ implies \
$ \,[\,   1 \!\perp \!\!\! \perp  3  \,] \,$ and
$ \,[\,  1 \!\perp \!\!\! \perp 2  \,\,|\, \{ 3   \} \,] \,$,
\item[(b)] $ \,[\,  1 \!\perp \!\!\! \perp  2 \,\,|\, \{  3  \} \,] \,$ and
$ \,[\,  1 \!\perp \!\!\! \perp 3  \,\,|\, \{ 2   \} \,] \,$ \ implies \
$ \,[\,  1 \!\perp \!\!\! \perp 2  \,] \,$ and
$ \,[\,   1 \!\perp \!\!\! \perp 3  \,] \,$,
\item[(c)] $ \,[\,   1 \!\perp \!\!\! \perp  2  \,] \,$ and
$ \,[\,  1 \!\perp \!\!\! \perp  3  \,] \,$  \ implies \
$ \,[\,   1 \!\perp \!\!\! \perp  2 \,\,|\, \{  3  \} \,] \,$ and
$ \,[\,  1 \!\perp \!\!\! \perp  3 \,\,|\, \{ 2  \} \,] \,$,
\item[(d)] $ \,[\, 1  \!\perp \!\!\! \perp 2   \,] \,$ and
$ \,[\,  1 \!\perp \!\!\! \perp  2 \,\,|\, \{ 3   \} \,] \,$ \ implies \
$ \,[\,  1 \!\perp \!\!\! \perp  3 \,] \,$ {\bf or}
$ \,[\,   2 \!\perp \!\!\! \perp 3  \,] \,$.
\end{itemize}
The necessity of these axioms can be
checked by simple calculations involving
almost-principal minors of  positive definite
symmetric $3 {\times} 3$-matrices:
\begin{itemize}
\item[(a)]  $\sigma_{12} = \sigma_{13} \sigma_{22} - \sigma_{12} \sigma_{23} = 0$ \ implies \
                  $\sigma_{13} = \sigma_{12} \sigma_{33} - \sigma_{13} \sigma_{23} = 0 $,
\item[(b)] $\sigma_{12} \sigma_{33} - \sigma_{13} \sigma_{23} = 
                \sigma_{13} \sigma_{22} - \sigma_{12} \sigma_{23} = 0 $ \ implies \
                 $\sigma_{12} = \sigma_{13} = 0$,
\item[(c)] $\sigma_{12} = \sigma_{13} = 0$ \ implies \
$\sigma_{12} \sigma_{33} - \sigma_{13} \sigma_{23} = 
                \sigma_{13} \sigma_{22} - \sigma_{12} \sigma_{23} = 0 $, 
\item[(d)] $\sigma_{12} = \sigma_{12} \sigma_{33} - \sigma_{13} \sigma_{23} = 0 $ \ implies \
                $\sigma_{13} = 0$ \ {\bf or} \ $\sigma_{23} = 0$.
\end{itemize}
The sufficiency of these axioms was noted in \cite[Example 1]{Mat}.

For arbitrary $n \geq 3$, a collection of
almost-principal minors is called a {\em gaussoid} if
it satisfies the axioms (a)-(d), after relabeling and applying
Schur complements. For instance, axiom (a) is then written as follows:
$ \,[\,  i \!\perp \!\!\! \perp  j  \,\,| \, L \,] $ and
$ [\,  i \!\perp \!\!\! \perp k  \,\,|\, \{ j   \} \cup L \,] $  implies 
$ [\,   i \!\perp \!\!\! \perp  k \,\,| \,L  \,] $ and
$ [\,  i \!\perp \!\!\! \perp j  \,\,|\, \{ k   \} \cup L \,] $.
This axiom is known as the {\em semigraphoid axiom}.
See \cite{MPSSW} for a discussion.

A gaussoid is {\em representable} if it is the set of
vanishing almost-principal minors of some matrix in ${\bf PD}_n$.
For $n=3$ every gaussoid is representable by \cite[Example 1]{Mat}.
For $n=4$, a complete classification of the representable gaussoids
was given in \cite{LM}. We are here asking for the extension to $n=5$.

We now introduce a  conceptual framework for our General Problem.
For each subset $S$ of $\{1,2,\ldots,n\}$ we introduce one unknown $H_S$,
and we define the {\em submodular cone} to be 
the solution set in $\RR^{2^n}$ of the system of
 linear inequalities
\begin{equation}
\label{submodIneq}
H_{\{i\} \cup K}+ H_{\{j\} \cup K} \, \,\leq\,\,
H_{\{i,j\} \cup K} \,+\,
H_{K},
\end{equation}
where $K$ is any subset of $\{1,\ldots,n\}$
and $i,j \in \{1,\ldots,n\} \backslash K$.
We denote this cone by $\,{\tt SubMod}_n \subset \RR^{2^n}$. Note that
 ${\tt SubMod}_n$ is a polyhedral cone living
in a high-dimensional space while
${\tt PD}_n$ is a non-polyhedral cone in a low-dimensional space.
Between these two cones we have the {\em entropy map}
$$ H :{\tt PD}_n \rightarrow {\tt SubMod}_n , $$
which is given by  the logarithms of all $\,2^n \,$ principal minors of 
 a positive definite matrix $\Sigma = (\sigma_{ij})$.
Namely, the coordinates of the entropy map are
 $$ H(\Sigma)_I \,\,\, = \,\,\, - {\rm log} \, {\rm det} (\Sigma_{I}), $$
 where $I $ is any subset of $\{1,\ldots,n\}$
 and $\Sigma_{I}$ the corresponding principal minor.
Note that the entropy map is well-defined because of the inequality
\begin{equation}
\label{DetIneq}
 {\rm det} (\Sigma_{\{i\} \cup K }) \cdot
{\rm det} (\Sigma_{\{j\} \cup K }) \,\,\, \geq \,\,\,
{\rm det} (\Sigma_{\{i,j\} \cup K }) \cdot
{\rm det} (\Sigma_{K }) .
\end{equation}
A matrix $\Sigma \in {\tt PD}_n$ satisfies (\ref{CIeqn})
 if and only if equality holds in (\ref{DetIneq}) 
 if and only if equality holds in (\ref{submodIneq}). 
This implies the following result.

\smallskip

\begin{proposition}
 The Gaussian conditional independence models  
are those subsets of the  positive definite cone $\,{\tt PD}_n\,$
that arise as inverse images of the  faces
of the submodular cone $\,{\tt SubMod}_n \,$
 under the entropy map $H$.
\end{proposition}

\smallskip

The importance of the submodular cone for probabilistic inference
with discrete random variables was highlighted in \cite{MPSSW}.
Here we are concerned with Gaussian random variables, and it
is the geometry of the entropy map which we must study.
We can thus paraphrase our problem as follows.

\smallskip

\noindent {\bf General Problem: } {\em Characterize the image of  the entropy map $\,H\,$
and how it intersects the various faces of ${\tt SubMod}_n$.
Study the fibers of this map.}

\smallskip

One approach to this problem is to work with the algebraic equations
satisfied by the principal minors of a symmetric matrix.
A characterization of these relations in terms of {\em hyperdeterminants}
was proposed in \cite{HS}. What we are interested in here
is the logarithmic image (or {\em amoeba}) of the positive
part of the hyperdeterminantal variety of \cite{HS}. A reasonable
first approximation to this amoeba is the tropicalization of that variety.
More precisely, we seek to compute the  {\em positive tropical variety}
\cite[\S 3.4]{ASCB} parametrically represented by the principal minors
of a symmetric $n {\times} n$-matrix.

\section{Bonus Problem on Rational Points} Section 4 dealt with conditional independence (CI) models
for Gaussians. Our bonus problem concerns CI models
for discrete random variables, thus returning to the setting
of Section 2. Consider $\,n \,$  discrete random variables
$X_1,X_2,\ldots,X_n\, $ with
$\,d_1,d_2,\ldots,d_n$ states. 
Any collection of CI statements $\, X_i \!\perp \!\!\! \perp X_j\,| X_K \,$
specifies a determinantal variety in the space of tables
\begin{equation}
\label{SpaceOfTables}
 \,\CC^{d_1} \otimes
 \CC^{d_2} \otimes \cdots \otimes
 \CC^{d_n} . 
 \end{equation}
 We call such a variety a {\em CI variety}. It is the zero set of
 a large collection of $2 {\times} 2$-determinants. These constraints are
 well-known and  listed explicitly in \cite[\S 4.1]{GMS} or \cite[Proposition 8.1]{Solving}.
The corresponding {\em strict CI variety} is the set 
of tables for which the given CI statements 
hold  but all other CI statements do not hold.
Thus a strict CI variety is a constructible subset of (\ref{SpaceOfTables})
which is Zariski open in a CI variety. The corresponding {\em strict CI model}
is the intersection of the strict CI variety with the positive orthant.
It consists of all positive $d_1 {\times} d_2 {\times} \cdots {\times} d_n$-tables
that lie in a common equivalence class, where two tables
are equivalent if precisely the same CI statements
$\, X_i \!\perp \!\!\! \perp X_j\,| X_K \,$ are valid (resp.~not valid) for both tables.

\smallskip

\noindent  {\bf Bonus Problem: } {\em Does every strict CI model have a $\QQ$-rational point?}

\smallskip

This charming problem was proposed by F.~Mat\'u\v{s} in \cite[page 275]{MatFinal}.
It suggests that algebraic statistics has something to offer also for
arithmetic geometers. One conceivable solution to the Bonus Problem might say
that CI models with no rational points exist but that 
rational points always appear when the number of states grows large,
that is, for  $d_1,d_2,\ldots, d_n \gg 0$. But that is pure speculation.
At present we know next to nothing.

\smallskip

\section{Brief Conclusion}
This article offered a whirlwind introduction to the emerging field
of algebraic statistics, by discussing a few of its numerous
open problems.  Aside from the Bonus Problem above, we had
listed three Specific Problems 
whose solution might be particularly rewarding:

\begin{itemize}
\item Consider the 
variety of $4 {\times} 4 {\times} 4$-tables
of tensor rank at most $4$.
Do the known polynomial invariants
of degree at most nine suffice
to define this variety? Set-theoretically? Ideal-theoretically? 
\item Characterize all projective varieties
whose   maximum likelihood degree is equal to one.
\item  Which sets of almost-principal minors can be simultaneously zero  
for a positive definite symmetric $5 {\times} 5$-matrix?
\end{itemize}

\end{document}